\documentclass[11pt,a4paper]{article}
\usepackage[english]{babel}
\usepackage{color}
\usepackage{amsfonts}
\usepackage{amssymb}
\usepackage{amsmath}
\usepackage{hyperref}
\usepackage[utf8]{inputenc}
\usepackage[T1]{fontenc}
\usepackage{graphicx}
\usepackage{url}
\newcommand{\bewende}{\hspace*{\fill}$\Box $\\ \vspace{0.5ex}}
\newcommand{\bewanf}{\noindent {\bf Proof:}\quad }

\DeclareMathOperator{\Li}{Li}
\newcommand{\R}{{\mathbb R}}
\newcommand{\N}{{\mathbb N}}

\newcommand{\C}{{\mathbb C}}
\newcommand{\satzanf}{\begin{samepage}\begin{satz}}
\newcommand{\satzende}{\end{satz}\end{samepage}}
\newcommand{\koranf}{\begin{samepage}\begin{kor}}
\newcommand{\korende}{\end{kor}\end{samepage}}
\newcommand{\lemanf}{\begin{samepage}\begin{lem}}
\newcommand{\lemende}{\end{lem}\end{samepage}}
\newcommand{\bspanf}{\begin{beispiel}}
\newcommand{\bspende}{\end{beispiel}}
\newcommand{\defanf}{\begin{definition}}
\newcommand{\defende}{\end{definition}}
\newtheorem{satz}{Theorem}
\newtheorem{lem}[satz]{Lemma}
\newtheorem{beispiel}[satz]{Example}
\newtheorem{definition}[satz]{Definition}
\newtheorem{kor}[satz]{Corollary}
\newcommand{\be}{\begin{equation}}
\newcommand{\ee}{\end{equation}}
\newcommand{\bl}[1]{\begin{equation}\label{#1}}
\newcommand{\bec}{\begin{center}}
\newcommand{\enc}{\end{center}}
\begin{document}

\author{Stephan Rosebrock}
\title{The Riemann Hypothesis and Prime Numbers\\ A simple Introduction}
\maketitle

\begin{abstract} 
The Riemann hypothesis is attracting a great deal of attention at present. It is one of the seven Millennium Prize Problems in mathematics, for each of which the Clay Mathematics Institute offers a prize of one million dollars. It is arguably the most important open problem in all of mathematics. Numerous results have so far been proved only under the assumption of the Riemann hypothesis. Yet the hypothesis itself sounds rather dull: it conjectures that the zeros of a certain function have a certain form. Things become truly exciting, however, when we look at its connection with prime numbers, and that is what we shall do here.

We give an elementary introduction to the subject.
\end{abstract}

\thispagestyle{empty}

\section{Prime Numbers}

We begin with a definition: a natural number $a$ {\it divides} a natural number $b$ if there is a natural number $c$ such that $b=a\cdot c$. We write $a|b$.

Prime numbers are natural numbers with exactly two positive divisors. Thus, the divisors of a prime number $p$ are 1 and $p$. It is known that there are infinitely many prime numbers.

The first prime numbers are:
\[2,3,5,7,11,13,17,19,23,29,31,37,41,43,47,53,\ldots\]
We shall need the following theorem, also known as the fundamental theorem of arithmetic.

\satzanf\label{eindeutprim}
Every natural number $n\geq 2$ can be expressed as a product of prime numbers. Apart from the order of the factors, this representation is unique: $n$ has {\em exactly one} prime factorization.
\satzende
For example, we can write $60=2^2\cdot3\cdot5$. In general, every natural number $n\geq2$ can be written as $n=p_1^{k_1}p_2^{k_2}\ldots p_j^{k_j}$, where the $p_i$ are distinct primes and $k_i\in\N$.

At first sight, the primes seem completely disordered. Sometimes they lie close together; at other times there are larger gaps. For example, one can find very large twin primes, that is, pairs of primes whose difference is 2. Whether there are infinitely many twin primes is still an open question. On the other hand, one can find arbitrarily long sequences of consecutive natural numbers, none of which is prime: The numbers $n!+2,\,n!+3,\ldots,n!+n$ are all nonprime for $n\ge2$ because $n!+k$ is divisible by $k$.

Here is a {\sc Julia} program for finding twin primes (further information on the {\sc Julia} programming language can be found in \cite{Jul26}).

\begin{verbatim}
# Simple primality test
function isprime(n::Int)
    if n < 2  return false end
    if n % 2 == 0  return n == 2  end
    limit = floor(Int, sqrt(n))
    for i in 3:2:limit
        if n % i == 0 return false end
    end
    return true
end

# Next pair of twin primes greater than n
function next_twin_prime(n::Int)
    p = max(n+1, 2)
    # Make p odd.
    if p % 2 == 0 p += 1 end
    while true
        if isprime(p) && isprime(p+2)
            return (p, p+2)
        end
        p += 2
    end
end
\end{verbatim}
The program is self-explanatory. Call it with any natural number, and it returns the next pair of twin primes above that number.

\begin{verbatim}
julia> next_twin_prime(3123425343452376402)
(3123425343452376749, 3123425343452376751)
\end{verbatim}

It therefore seems as though the primes have no structure. Tables can give us an impression of their distribution. Let $\Pi(n,m)$ denote the number of primes between $n$ and $m$.
\bec
\begin{tabular}{c|c}
$n - m$ & $\Pi(n,m)$\\
\hline
0 -- 100000 & 9592\\
100000 -- 200000 & 8392\\
200000 -- 300000 & 8013\\
300000 -- 400000 & 7863\\
400000 -- 500000 & 7678\\
500000 -- 600000 & 7560\\
600000 -- 700000 & 7445\\
700000 -- 800000 & 7408\\
800000 -- 900000 & 7323\\
900000 -- 1000000 & 7224\\
\end{tabular}
\enc

This suggests that there may be some underlying structure after all. The average density decreases in the long run. How exactly does it decrease?\\

In 1792, at the age of fifteen, {\sc Carl Friedrich Gauss} studied tables of primes, which could then be purchased as booklets, to understand these distributions. He formulated a conjecture that was proved only a hundred years later by Hadamard and, independently, by de la Vall\'{e}e Poussin. This conjecture became known as the ``prime number theorem.'' In what follows, let $\pi(n)$ denote the number of primes less than or equal to $n$. For example, $\pi(5)=\pi(6)=3$, since there are three primes less than or equal to either 5 or 6.

The following {\sc Julia} program, \verb+prime_pi(x)+, computes the function $\pi(x)$. It uses the function \verb+pi_tabelle(N)+, which applies the sieve of Eratosthenes to build a table of prime counts up to a given number \verb+N+.

\begin{verbatim}
function pi_tabelle(N::Integer)
    is_prime = trues(N)
    is_prime[1] = false  # 1 is not prime.

    for p in 2:isqrt(N)
        if is_prime[p]
            for multiple in (p * p):p:N
                is_prime[multiple] = false
            end
        end
    end

    return cumsum(vcat(0, Int.(is_prime)))
end
\end{verbatim}

\pagebreak

\begin{verbatim}
"""Number of primes less than or equal to a real number x."""
function prime_pi(x::Real)
    x < 2 && return 0
    return pi_tabelle(floor(Int, x))[end]
end
\end{verbatim}
The central line of the sieve is
\begin{verbatim}
for multiple in (p * p):p:N   is_prime[multiple] = false
\end{verbatim}
Starting at $p\cdot p$ and continuing up to $N$, every $p$th entry is crossed out, that is, set to \verb+false+ in \verb+is_prime+.

\begin{figure}[htbp]
    \centering
    \includegraphics[width=0.7\textwidth]{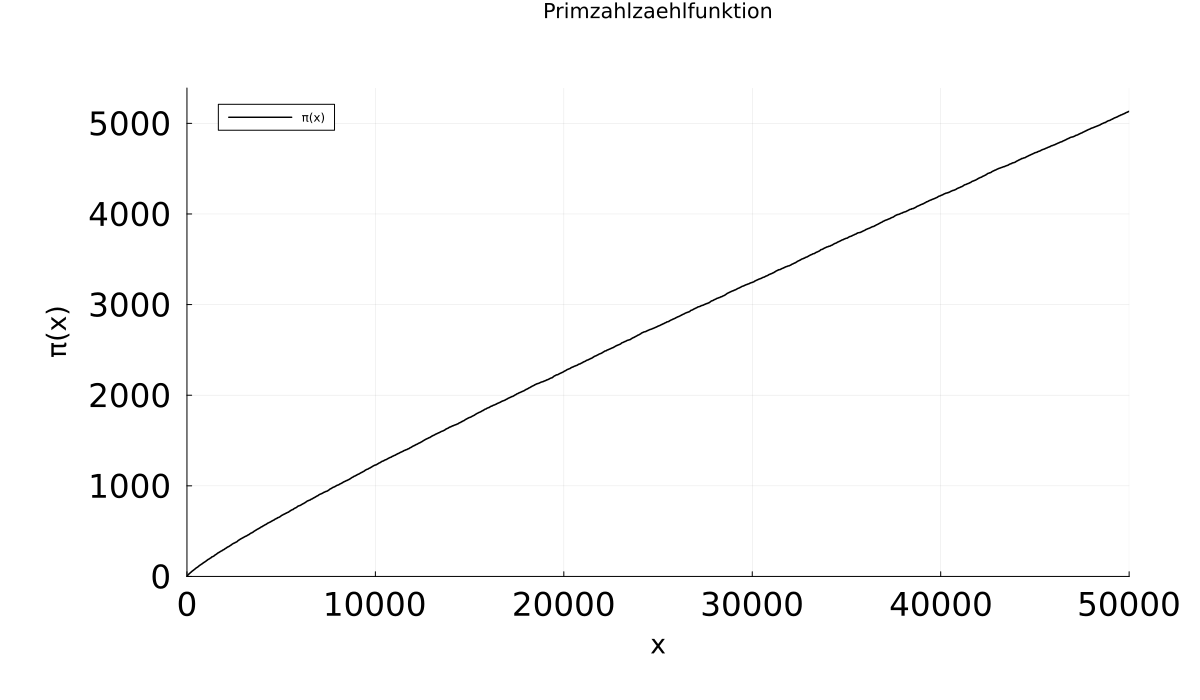}
    \caption{The prime-counting function $\pi(x)$ for $0\leq x\leq50\,000$.}
    \label{fig:primzahlfunktion}
\end{figure}

Plotting the graph of $\pi(x)$ gives Figure~\ref{fig:primzahlfunktion}. The regularity of the curve is astonishing. There must be some structure behind it.

The following table lists the values of $x/\pi(x)$. We observe that, for larger values of $x$, they increase by about $2.3$ from one power of ten to the next.

\begin{tabular}{r|r|r}
$x$ & $\pi(x)$ & $x/\pi(x)$ \\
\hline
10 & 4 & $2.5$ \\
100 & 25 & $4.0$ \\
1000 & 168 & $6.0$ \\
10000 & 1229 & $8.1$ \\
100000 & 9592 & $10.4$ \\
1000000 & 78498 & $12.7$ \\
10000000 & 664579 & $15.0$ \\
100000000 & 5761455 & $17.4$ \\
1000000000 & 50847534 & $19.7$ \\
10000000000 & 455052511 & $22.0$ \\
\end{tabular}

\pagebreak

We have $\ln(10)\approx2.3$, where $\ln x$ denotes the natural logarithm of $x$. This suggests:
\[\frac{10^n}{\pi(10^n)}\approx n\cdot\ln(10)=\ln(10^n).\]
If this holds for $10^n$, perhaps it holds more generally:
$x/\pi(x)\approx\ln(x)$, and hence $\pi(x)\approx x/\ln(x)$. This was Gauss's discovery.

\satzanf\label{sgrPr} {\bf Prime Number Theorem}
As $x$ increases, $\pi(x)$ approaches the function $f$ given by $f(x):=\frac{x}{\ln x}$ in the following sense:
$$\lim_{x\to\infty}\frac{\pi(x)}{{x\over\ln x}}=1.$$
\satzende

One can interpret this as follows: the ``probability'' that a randomly chosen number $n$ is prime is approximately
\[\frac{1}{\ln(n)}\approx\frac{1}{2.3\cdot\mbox{number of digits}}.\]

The relative frequency of primes among all numbers up to $x$, namely $\pi(x)/x$, is therefore approximately $1/\ln(x)$. This leads to another good approximation to $\pi(x)$, first proposed by Gauss. The idea is to assign each number a weight corresponding to the approximate density of primes in its neighborhood.

First, we can write the prime-counting function exactly as a sum. For natural numbers $x\ge2$,
\bl{gan}
\pi(x)=\sum_{n=2}^{x}a_n,
\qquad
a_n=
\begin{cases}
1,&\text{if }n\text{ is prime},\\
0,&\text{otherwise}.
\end{cases}
\ee

The individual values $a_n$ jump irregularly between $0$ and $1$. Gauss's idea is to replace these irregular contributions by smooth average contributions. His heuristic is:
\begin{quote}
Near a large number $n$, approximately a proportion $1/\ln n$ of the numbers are prime.
\end{quote}

For example,
\[\frac{1}{\ln1000}\approx0.145.\]
Intuitively, this means that in a suitable range near $1000$, we expect about $14.5\,\%$ of the numbers to be prime. For a fixed number, of course, whether it is prime is completely determined; the ``probability'' here describes an average frequency.

We therefore approximately replace the contribution $a_n$ in (\ref{gan}) by $1/\ln n$. Adding these weights gives
\[
\pi(x)\approx\frac{1}{\ln2}+\frac{1}{\ln3}+\cdots+\frac{1}{\ln x}.
\]

This takes into account the fact that prime density decreases as $n$ increases: smaller numbers receive larger weights, and larger numbers receive smaller weights. The simpler approximation $x/\ln x$ uses the same value $1/\ln x$ throughout the whole range.

\begin{figure}[htbp]
    \centering
    \includegraphics[width=0.8\textwidth]{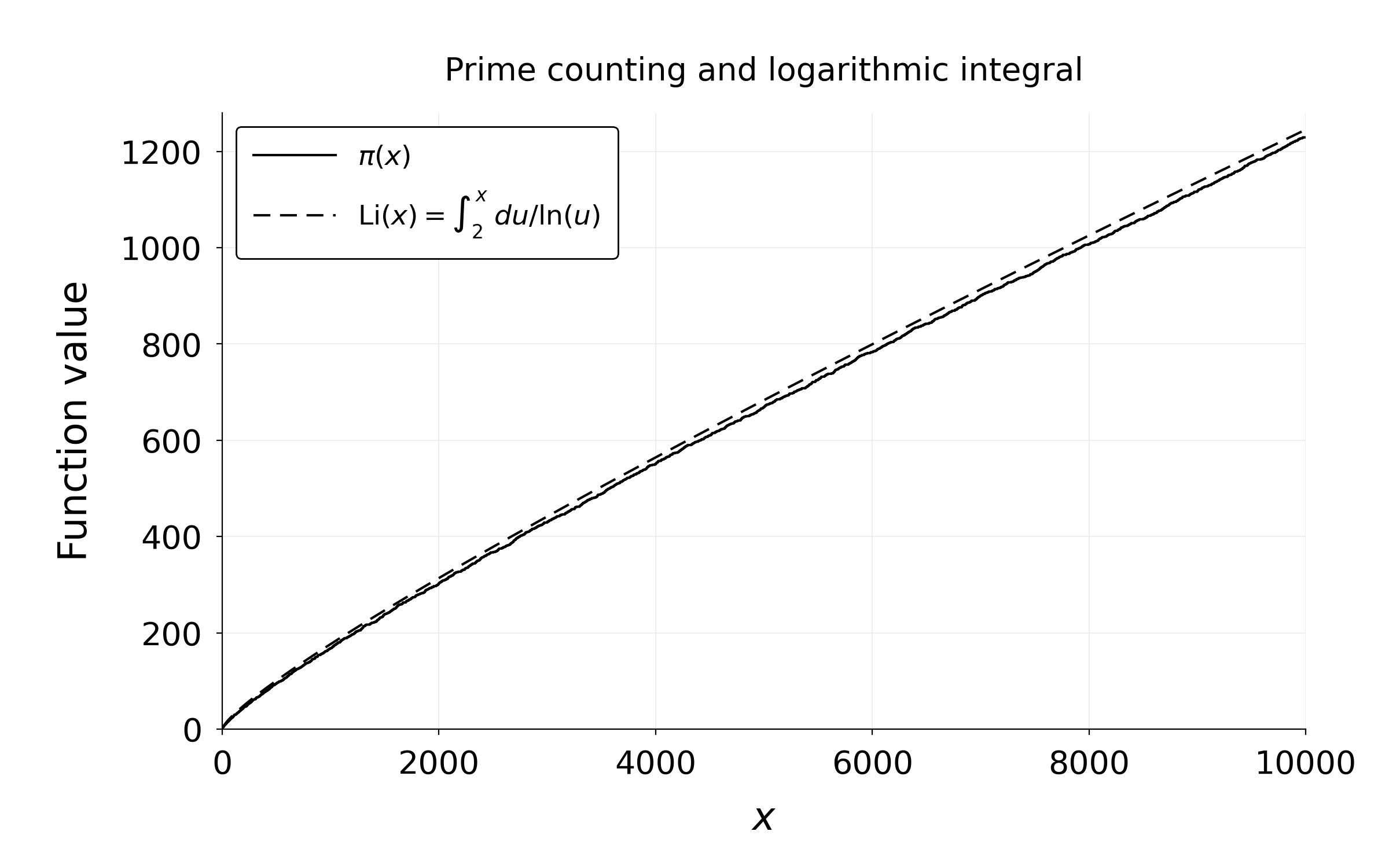}
    \caption{$\pi(x)$ and $\Li(x)$ for $2\leq x\leq10\,000$.}
    \label{fig:Lipi}
\end{figure}

The sum can in turn be approximated by the area under the graph of $u\mapsto1/\ln u$:
\[
\sum_{n=2}^{x}\frac{1}{\ln n}\approx\int_2^x\frac{du}{\ln u}.
\]
We therefore expect
$$\pi(x)\approx\int_2^x\frac{du}{\ln u}.$$
This approximation is somewhat more accurate than the one in the prime number theorem. We define
\[\Li(x)=\int_2^x\frac{du}{\ln u}.\]
Plotting $\Li$ and $\pi$ together gives Figure~\ref{fig:Lipi}.
Notice that $\pi(x)$ always lies below $\Li(x)$ (at least for $x>8$). This remains the case when the graph is extended to much larger values of $x$. However, Littlewood has proved in 1914 that $\pi(x)$ and $\Li(x)$ exchange places infinitely often. The first such crossing occurs for $x<10^{316}$. We may never find this point: $10^{300}$ is simply too large.

The following statement is equivalent to the Riemann hypothesis (be patient---we shall introduce it later):

\begin{quote}
{\bf (RH'):} For sufficiently large $x$, the absolute error in Gauss's approximation $\Li(x)$ to $\pi(x)$ is at most a constant multiple of $\sqrt{x}\cdot\ln x$, that is,
\[
|\pi(x)-\Li(x)|\le C\sqrt{x}\ln x.
\]
\end{quote}
This can also be made more explicit. For example, the Riemann hypothesis is equivalent to
\[
|\pi(x)-\Li(x)|<\frac{1}{8\pi}\sqrt{x}\ln x,\qquad\text{for }x\ge2657.
\]
Schoenfeld proved this as early as 1976. This remains an active area of research.

\section{A Formula of Euler}

Let $1\ne q^s$. It is easy to rove the formula
\bl{gf1}
\frac{1-q^{(n+1)s}}{1-q^s}=1+q^s+q^{2s}+q^{3s}+\ldots+q^{ns}.
\ee
Expand $(1-p^s)\cdot(1+p^s+p^{2s}+p^{3s}+\ldots+p^{ns})$ by hand and show that all terms cancel except $1$ and $-p^{(n+1)s}$.

\lemanf\label{lunendS} For $p\in\N$, $p>1$, and $s>1$, we have
\[\sum_{i=0}^\infty\frac{1}{(p^s)^i}=\frac{1}{1-\displaystyle\frac{1}{p^s}}.\]
\lemende

\bewanf Substitute $q=1/p$ into formula (\ref{gf1}):
\[\sum_{i=0}^n\frac{1}{(p^s)^i}=\frac{1-\displaystyle\frac{1}{p^{(n+1)s}}}{1-\displaystyle\frac{1}{p^s}}.\]
As $n$ tends to infinity, $\frac{1}{p^{(n+1)s}}$ tends to 0 (since $p>1$), giving the claimed formula.
\bewende

Let $\cal P$ be the set of prime numbers. Euler considered the product
\[
\prod_{p\in\cal P}\frac{1}{1-\frac{1}{p^s}}=
\frac{1}{1-\frac{1}{2^s}}\cdot\frac{1}{1-\frac{1}{3^s}}\cdot\frac{1}{1-\frac{1}{5^s}}\cdot\ldots.
\]
By Lemma~\ref{lunendS}, we can rewrite this as
\[\prod_{p\in\cal P}\frac{1}{1-\frac{1}{p^s}}=\prod_{p\in\cal P}\sum_{i=0}^\infty\frac{1}{(p^s)^i}.\]

We obtain
\begin{align*}
\prod_{p\in\cal P}\frac{1}{1-\frac{1}{p^s}}
&=\prod_{p\in\cal P}\sum_{i=0}^\infty\frac{1}{(p^s)^i}\\
&=\left(1+\frac{1}{2^s}+\frac{1}{2^{2s}}+\frac{1}{2^{3s}}+\frac{1}{2^{4s}}+\ldots\right)\\
&\quad\cdot\left(1+\frac{1}{3^s}+\frac{1}{3^{2s}}+\frac{1}{3^{3s}}+\frac{1}{3^{4s}}+\ldots\right)\\
&\quad\cdot\left(1+\frac{1}{5^s}+\frac{1}{5^{2s}}+\frac{1}{5^{3s}}+\frac{1}{5^{4s}}+\ldots\right)\\
&\quad\cdot\left(1+\frac{1}{7^s}+\frac{1}{7^{2s}}+\frac{1}{7^{3s}}+\frac{1}{7^{4s}}+\ldots\right)\\
&\quad\cdot\ldots.
\end{align*}
In this product, the nonzero contributions are obtained by choosing a value different from 1 from only finitely many factors. Each such term has the form
\[
\frac{1}{p_1^{sk_1}p_2^{sk_2}\ldots p_j^{sk_j}}=
\left(\frac{1}{p_1^{k_1}p_2^{k_2}\ldots p_j^{k_j}}\right)^s,
\]
where the $p_i$ are distinct primes and $k_i\in\N$. Moreover, each possible term occurs exactly once.
Recall Theorem~\ref{eindeutprim}, the theorem on unique prime factorization. The denominator on the right is the prime factorization of a natural number. Every natural number appears exactly once as a denominator. Thus we obtain the following formula, which Euler discovered in the eighteenth century:
\bl{gf2}
\prod_{p\in\cal P}\frac{1}{1-\frac{1}{p^s}}=\sum_{n=1}^\infty\frac{1}{n^s}.
\ee

This formula is another marvel. On the left is a product involving only primes. On the right, we sum over all natural numbers. At first sight, the expression on the right has nothing to do with primes at all! Yet somehow all the primes must be hidden in that sum. This must have caught Riemann's attention.

\section{The Riemann Zeta Function}

Formula (\ref{gf2}) defines a function of $s$. The function
\[\zeta(s)=\sum_{n=1}^\infty\frac{1}{n^s}\]
is called the {\it Riemann zeta function}. It is defined only for $s>1$. For $s=1$, we would have
\[\zeta(1)=\sum_{n=1}^\infty\frac{1}{n},\]
and here the following result applies:

\satzanf\label{srei} The series $\sum_{n=1}^\infty\frac{1}{n}$ diverges.\satzende

\bewanf We group the terms into blocks:
\[
1+\frac12
+\left(\frac13+\frac14\right)
+\left(\frac15+\cdots+\frac18\right)
+\left(\frac19+\cdots+\frac1{16}\right)+\cdots.
\]
The block
\[\sum_{n=2^{k-1}+1}^{2^k}\frac1n\]
contains $2^{k-1}$ terms. Each is at least $1/2^k$. Therefore,
\[
\sum_{n=2^{k-1}+1}^{2^k}\frac1n\ge\frac{2^{k-1}}{2^k}=\frac12.
\]
Hence the partial sums satisfy
\[
\sum_{n=1}^{2^m}\frac1n\ge1+\frac m2\xrightarrow[m\to\infty]{}\infty.
\]
Since all terms are positive, the partial sums increase monotonically. Thus
\[{\displaystyle\sum_{n=1}^{\infty}\frac1n=\infty.}\]
\bewende

Incidentally, this gives us a proof that there are infinitely many primes. If there were only finitely many, the product on the left of formula (\ref{gf2}) would be a finite number, whereas the right-hand side is infinite by the theorem just proved. This is a contradiction.

In 1859, Riemann wrote a groundbreaking paper \cite{Rie59}. His idea was to extend the function $\zeta(s)$ to complex numbers (the necessary background on complex numbers can be found, for example, in \cite{EnFe16}). This makes $\zeta(s)$ a function defined for all complex numbers $s$ whose real part is greater than 1.

We now have natural numbers raised to complex powers. How should we think about this?

Let $n\ge1$ be a natural number and let
\[s=\sigma+it,\qquad\sigma>1.\]
We define $n^s:=e^{s\ln n}$. Euler's formula (i.e. $e^{i\theta}=\cos \theta+i\sin \theta$ for $\theta\in \R$) gives
\[
n^{\sigma+it}=n^\sigma\bigl(\cos(t\ln n)+i\sin(t\ln n)\bigr).
\]
Thus the modulus is $n^\sigma$, while the angle is $t\ln n$.

For example, $2^{2+i}=4\bigl(\cos(\ln2)+i\sin(\ln2)\bigr)$.

For the reciprocal, we similarly have
\[
\frac1{n^s}=n^{-\sigma}\bigl(\cos(t\ln n)-i\sin(t\ln n)\bigr),
\qquad
\left|\frac1{n^s}\right|=\frac1{n^\sigma}.
\]
These powers are defined for every complex exponent. The condition $\sigma>1$ additionally guarantees that the series
\[\sum_{n=1}^{\infty}\frac1{n^s}\]
converges absolutely, since
\[
\sum_{n=1}^{\infty}\left|\frac1{n^s}\right|=\sum_{n=1}^{\infty}\frac1{n^\sigma}<\infty.
\]
\begin{figure}[htbp]
    \centering
    \includegraphics[width=0.6\textwidth]{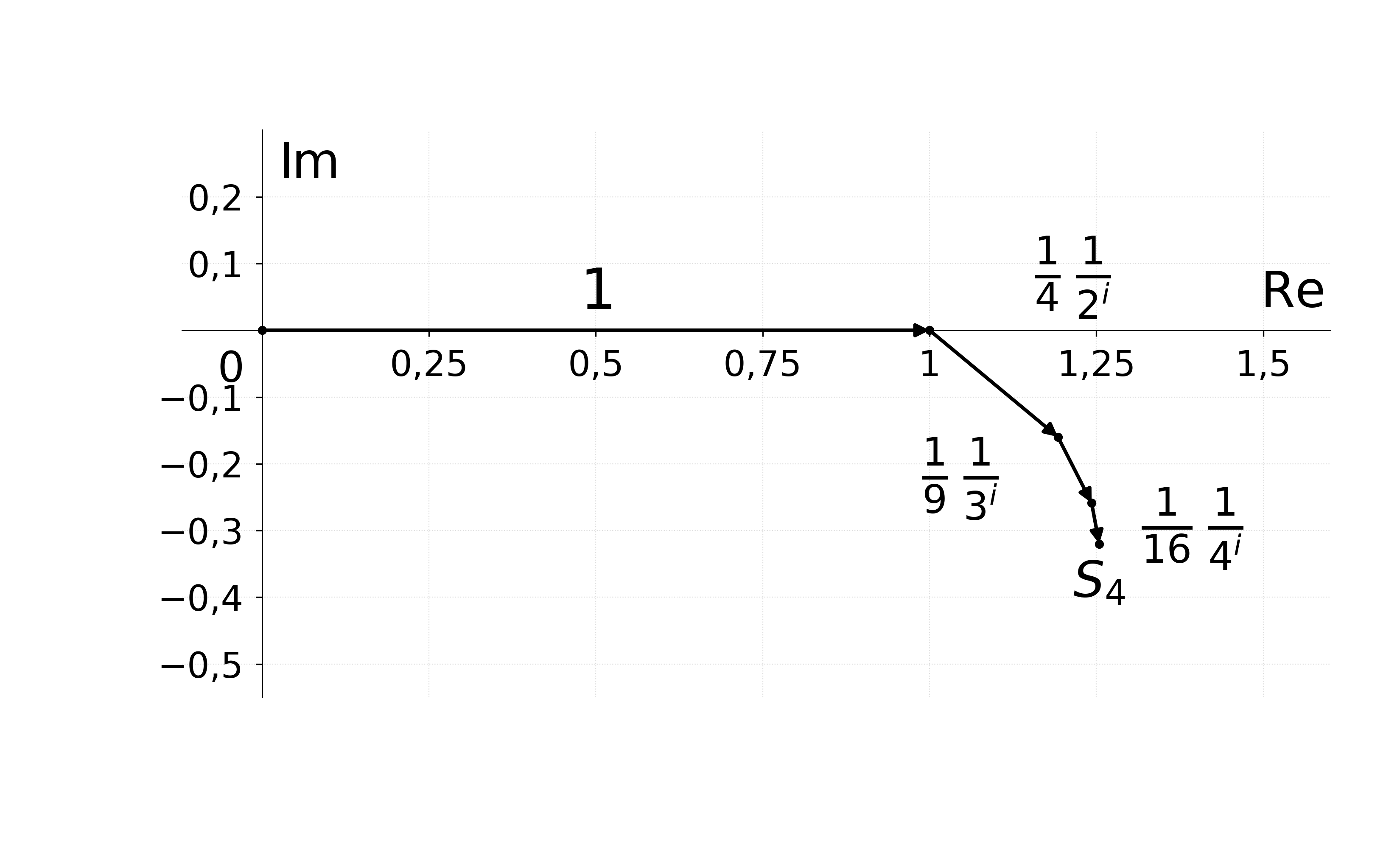}
    \caption{An approximation to $\zeta(2+i)$.}
    \label{fig:z2i}
\end{figure}
Let us look at a concrete example. Consider
\begin{align*}
\zeta(2+i)&=\frac{1}{1^{2+i}}+\frac{1}{2^{2+i}}+\frac{1}{3^{2+i}}+\frac{1}{4^{2+i}}+\dots\\
&=1+\frac14\frac{1}{2^i}+\frac19\frac{1}{3^i}+\frac1{16}\frac{1}{4^i}+\ldots.
\end{align*}
Figure~\ref{fig:z2i} shows an approximation to this function value in the complex plane. Taking only the first term would bring us to 1, as indicated by the arrow from 0 to 1. Adding the second term brings us to the tip of the second arrow. We have turned through $-39.7^\circ$, due to the factor $1/2^i$, and traveled a distance of $1/4$. The other two terms are handled in the same way. The point $S_4$ represents the function value obtained from the first four terms.

To summarize, we have a function $\zeta(s)=\sum_{n=1}^\infty\frac{1}{n^s}$ whose domain consists of the complex numbers with real part greater than 1. It is known---and Riemann already knew this---that $\zeta(s)$ has a unique holomorphic continuation (that is, a continuation that is complex differentiable everywhere) to the entire complex plane except for a simple pole at $s=1$.

\begin{figure}[htbp]
    \centering
    \includegraphics[width=0.6\textwidth]{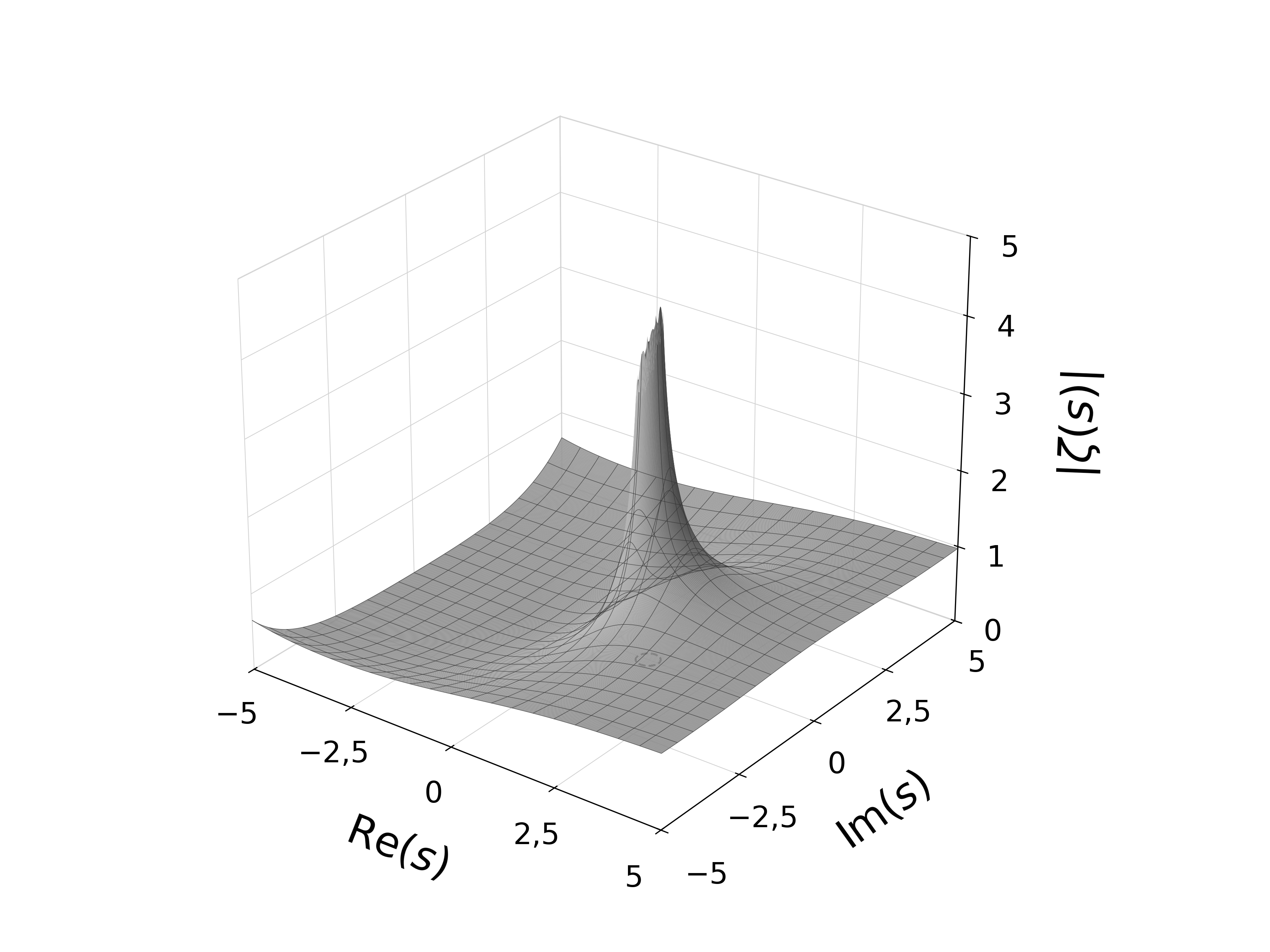}
    \caption{The modulus of the zeta function, with the pole at 1 excluded.}
    \label{fig:bz}
\end{figure}
Figure~\ref{fig:bz} shows the modulus of the zeta function extended to $\C\backslash\{1\}$. The pole is at 1, and we can see how the modulus grows in its neighborhood.

Hasse's series provides a formula for computing the zeta function throughout $\C\backslash\{1\}$ (see \cite{Has30}):
\[
\zeta(s)=
\frac{1}{s-1}
\sum_{n=0}^{\infty}\frac{1}{n+1}
\sum_{k=0}^{n}
(-1)^k\binom{n}{k}(k+1)^{1-s},
\qquad s\in\mathbb{C}\setminus\{1\}.
\]

Let us experiment with {\sc Julia}. The zeta function is implemented in Julia as \texttt{zeta(s)}. It belongs to the \texttt{SpecialFunctions} package and accepts both real and complex arguments.

Install the package once in Julia:
\begin{verbatim}
import Pkg; Pkg.add("SpecialFunctions")
\end{verbatim}
Then load it using
\begin{verbatim}
using SpecialFunctions
\end{verbatim}
You can now compute function values:
\begin{verbatim}
zeta(2)       # Approximately 1.6449340668482264
zeta(0)       # -0.5
zeta(-1)      # Approximately -0.08333333333333333
zeta(-2)      # 0.0
\end{verbatim}
In Julia, \texttt{im} denotes the imaginary unit $i$.
\begin{verbatim}
s = 0.5 + 14.1im
w = zeta(s)

println("Function value: ", w)
println("Real part:      ", real(w))
println("Imaginary part: ", imag(w))
println("Modulus:        ", abs(w))
\end{verbatim}

\section{The Riemann Hypothesis}

Consider the unique holomorphic continuation to $\mathbb C\setminus\{1\}$ of the function initially defined by $\sum_{n\geq1}n^{-s}$ for $\operatorname{Re}s>1$ (Re $s$ means the real part of the complex number $s$). We continue to denote it by $\zeta$.

The $\zeta$-function has no zeros with real part greater than 1. This can be seen from the fact that $\zeta(s)=\sum_{n=1}^\infty\frac{1}{n^s}$ is nonzero for real numbers $s>1$.
One can prove that the negative even integers, that is, $-2,-4,-6,\ldots$, are zeros of $\zeta$. They are called the {\it trivial zeros}. There are also nontrivial zeros; indeed, it is known that there are infinitely many of them. An online database contains more than 103 billion nontrivial zeros of the zeta function. All these zeros have real part $1/2$.

In Table \ref{t1} are the first ten nontrivial zeros ($\operatorname{Re}$ denotes the real part and $\operatorname{Im}$ the imaginary part of a complex number).
\begin{table}[htbp]
    \centering
    \renewcommand{\arraystretch}{1.15}
    \begin{tabular}{c|c|r}
        \hline
        $n$ & $\operatorname{Re}(\rho_n)$ & $\operatorname{Im}(\rho_n)$ \\
        \hline
         1 & $1/2$ & $14.1347251417$ \\
         2 & $1/2$ & $21.0220396388$ \\
         3 & $1/2$ & $25.0108575801$ \\
         4 & $1/2$ & $30.4248761259$ \\
         5 & $1/2$ & $32.9350615877$ \\
         6 & $1/2$ & $37.5861781588$ \\
         7 & $1/2$ & $40.9187190121$ \\
         8 & $1/2$ & $43.3270732809$ \\
         9 & $1/2$ & $48.0051508812$ \\
        10 & $1/2$ & $49.7738324777$ \\
        \hline
    \end{tabular}
    \caption{\label{t1}The first ten nontrivial zeros $\rho_n=\frac12+i\gamma_n$ of the Riemann zeta function with positive imaginary part. The imaginary parts are rounded to ten decimal places.}
    \label{tab:zeta-nullstellen}
\end{table}
For every nontrivial zero, its complex conjugate, with the opposite imaginary part, is also a zero. For example, $1/2-14.1347251417i$ is also a zero.

We can now state the Riemann hypothesis:
\begin{quote}
{\bf (RH):} All nontrivial zeros of the Riemann zeta function have real part $1/2$.
\end{quote}
In the first section, we encountered conjecture (RH'), a statement about primes. Conjecture (RH) is an equivalent statement in analysis. The zeta function translates arithmetic statements about primes into analytic properties of functions. The connection is astonishing. In modern mathematics, one often tries to reformulate conjectures from one field in another in order to obtain a different approach, as in the Langlands program -- but that is another story.

The line $x=1/2$ in the complex plane is called the {\it critical line}. The Riemann hypothesis can therefore also be stated as follows:

{\it All nontrivial zeros of the Riemann zeta function lie on the critical line.}

It is known that all nontrivial zeros lie in the {\it critical strip}, that is, $0<\operatorname{Re}s<1$ for every nontrivial zero $s$ (these bounds have recently been slightly improved).
\begin{figure}[htbp]
    \centering
    \includegraphics[width=0.6\textwidth]{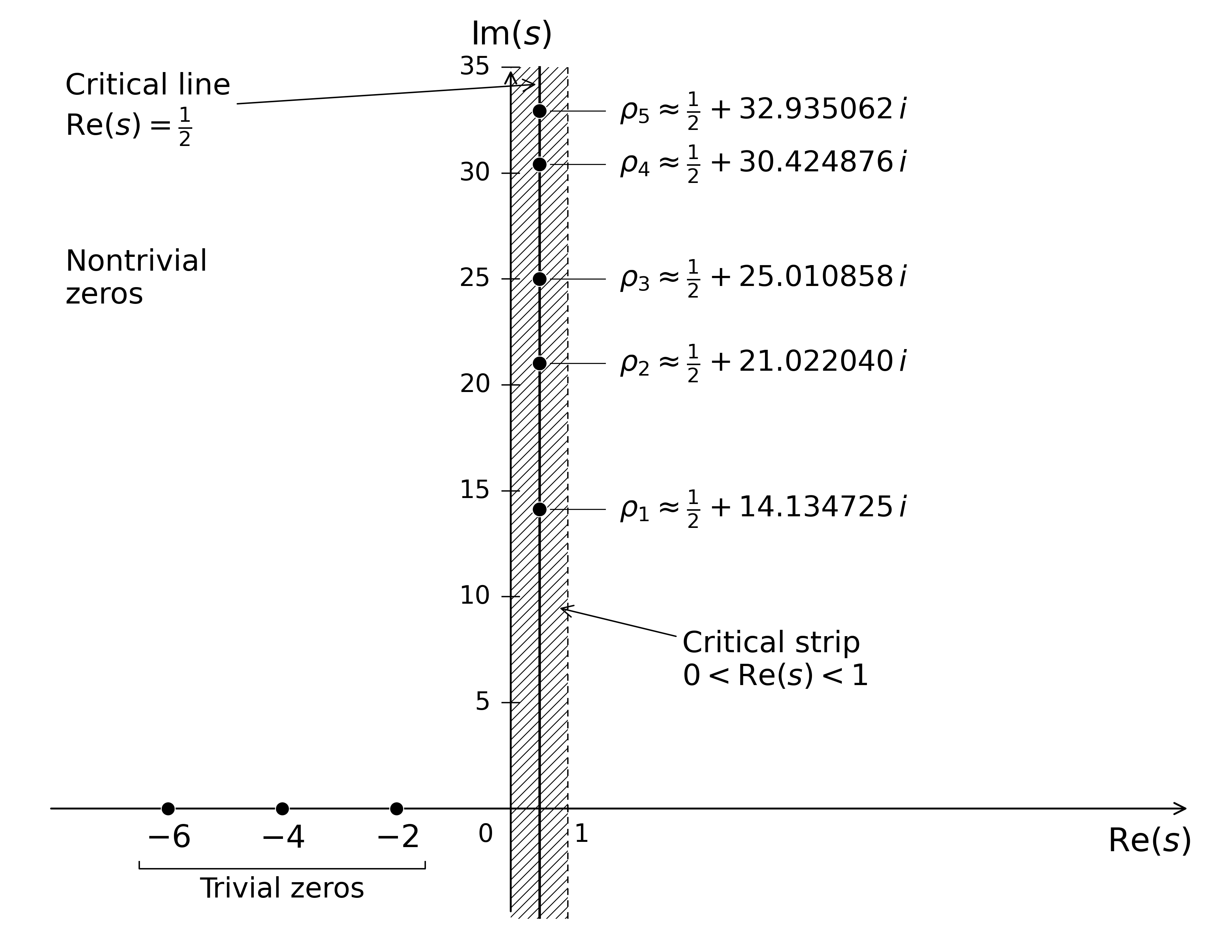}
    \caption{The complex plane.}
    \label{fig:cE}
\end{figure}
Figure~\ref{fig:cE} shows part of the complex plane with the first trivial and nontrivial zeros.

\begin{figure}[htbp]
    \centering
    \includegraphics[width=0.7\textwidth]{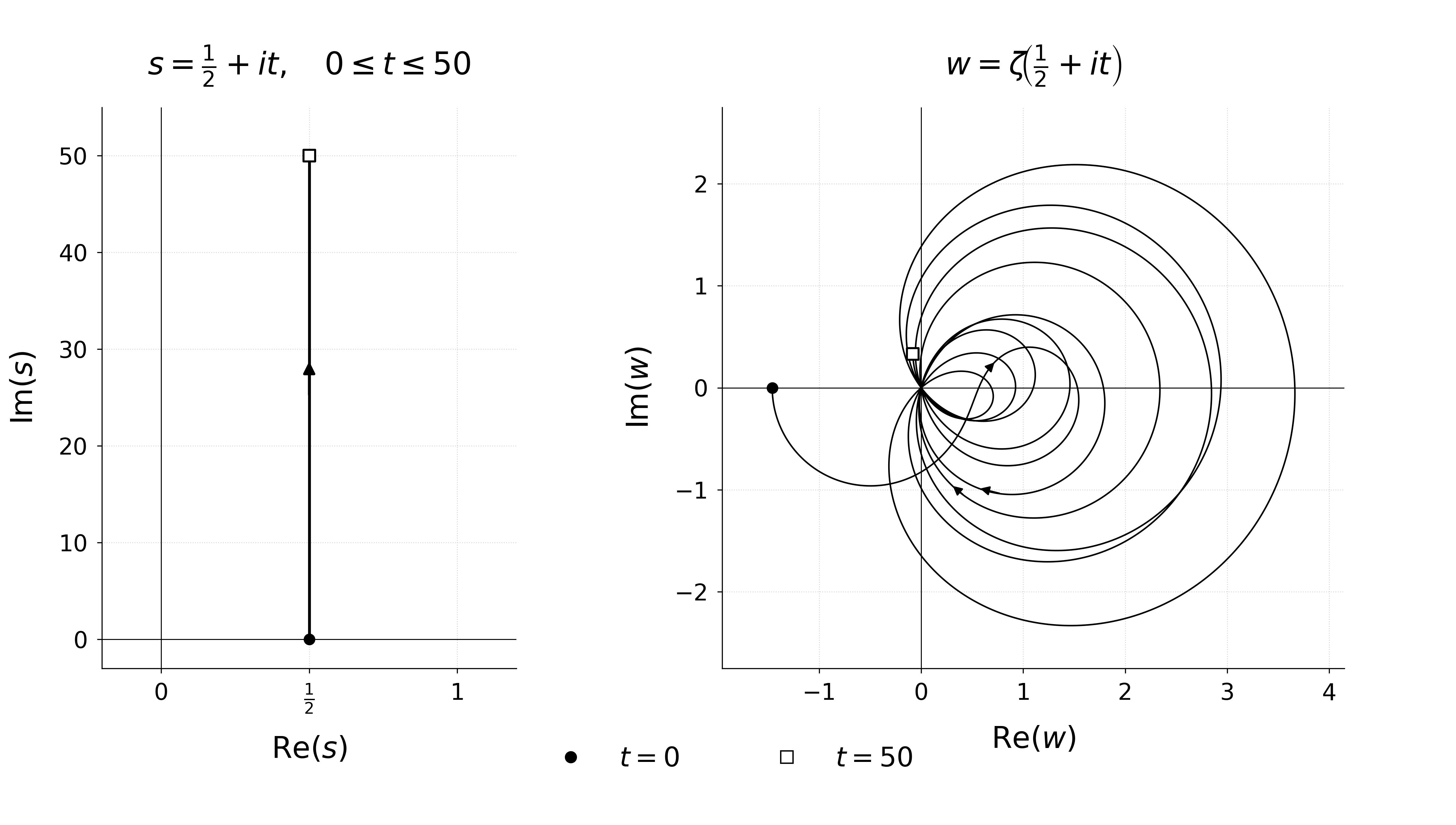}
    \caption{The segment $s=1/2+it$, $0\le t\le50$, in the domain and its image under $\zeta$.}
    \label{fig:zkg}
\end{figure}
Figure~\ref{fig:zkg} shows a segment of the critical line and its image under $\zeta$. We can see the image curve repeatedly passing through 0.

The Riemann hypothesis has many consequences, for example for primality tests, elliptic curves, endomorphism rings, and much more. In his paper \cite{Rie59}, Riemann himself developed an analytic formula for $\pi(x)$. In the next section, we turn to another application.

\section{Wave Functions}

The followiung can be found in \cite[Chapters~30, 32, and~35]{MaSt16}.
Let $\gamma_i$ be the imaginary parts of the nontrivial zeros of the zeta function in the upper half-plane, arranged in increasing order. Thus $\gamma_1=14.1347251417$, $\gamma_2=21.0220396388$, and so on. Consider the functions $f_i(x)=\cos(\gamma_i\cdot\ln(x))$.

\begin{figure}[htbp]
    \centering
    \includegraphics[width=0.6\textwidth]{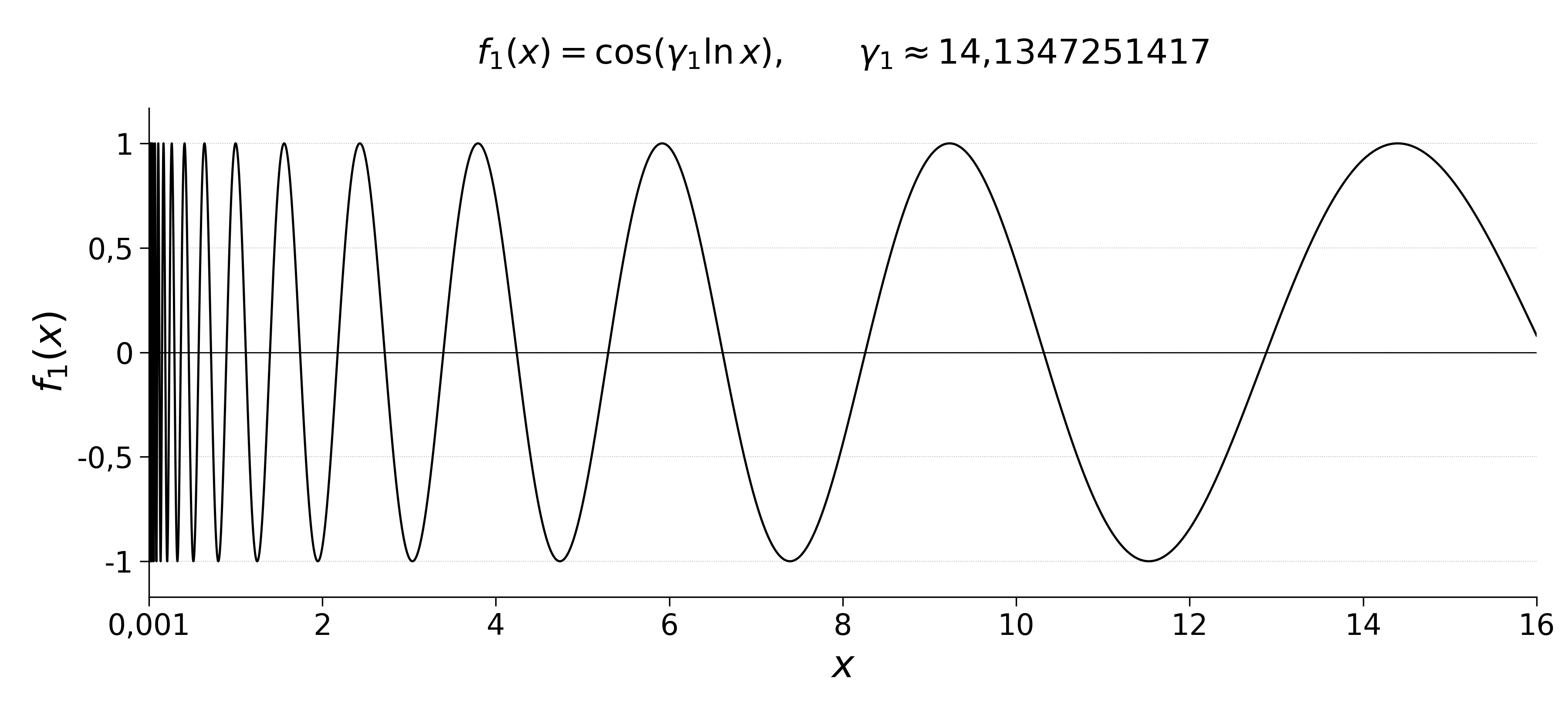}
    \caption{The function $f_1(x)=\cos(\gamma_1\cdot\ln(x))$.}
    \label{fig:f1}
\end{figure}

Figure~\ref{fig:f1} shows the function $f_1(x)=\cos(\gamma_1\cdot\ln(x))=\cos(14.1347251417\cdot\ln(x))$.

\begin{figure}[htbp]
    \centering
    \includegraphics[width=0.6\textwidth]{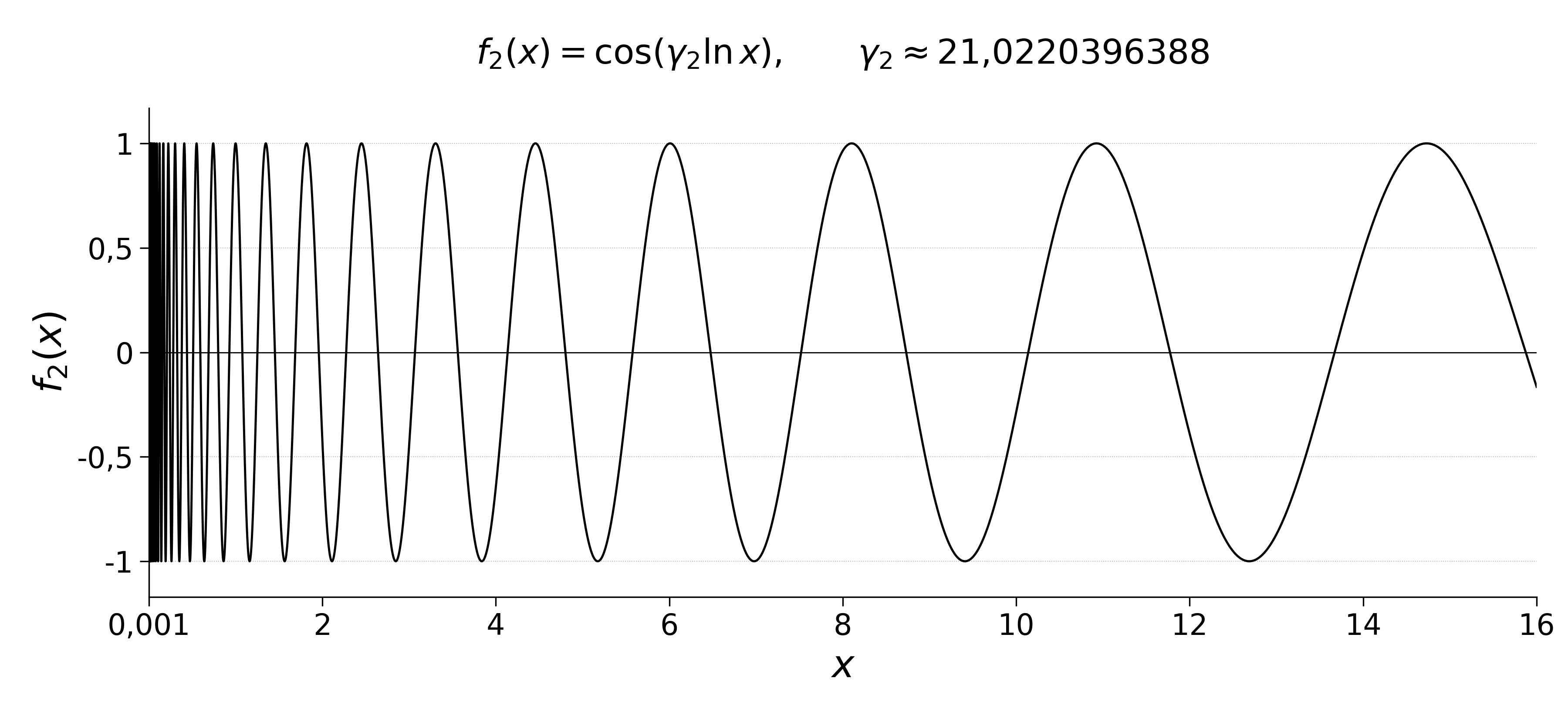}
    \caption{The function $f_2(x)=\cos(\gamma_2\cdot\ln(x))$.}
    \label{fig:f2}
\end{figure}

\begin{figure}[htbp]
    \centering
    \includegraphics[width=0.9\textwidth]{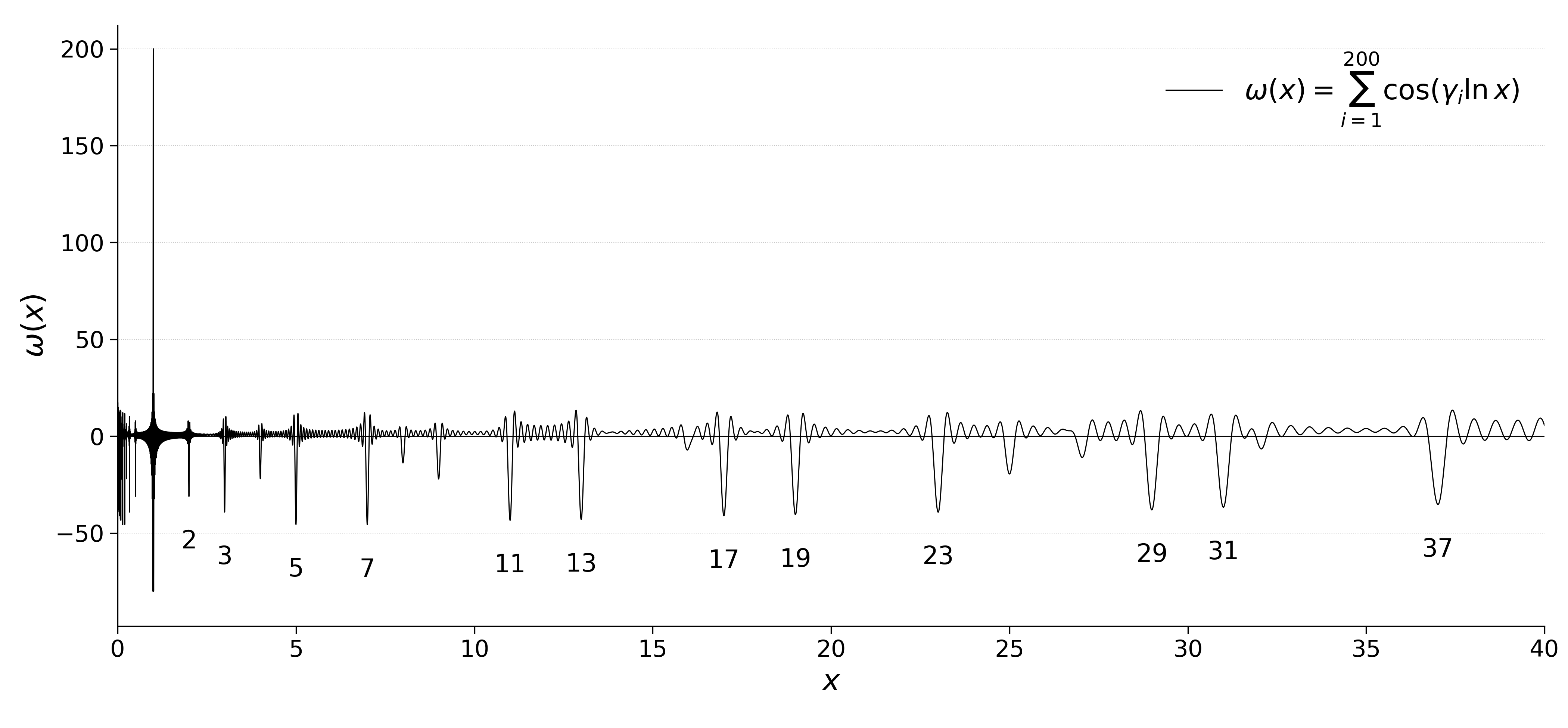}
    \caption{$f_1(x)+f_2(x)+\ldots+f_{200}(x)$.}
    \label{fig:cossum}
\end{figure}

Figure~\ref{fig:f2} shows the function $f_2(x)=\cos(\gamma_2\cdot\ln(x))=\cos(21.0220396388\cdot\ln(x))$.
It is very similar to $f_1$, but somewhat compressed horizontally. We can superimpose these wave functions by adding them:
\[\omega(x)=\sum_{i=1}^{200}f_i(x)=f_1(x)+f_2(x)+\ldots+f_{200}(x).\]
The resulting function $\omega(x)$ is plotted in Figure~\ref{fig:cossum}.
There are pronounced downward spikes very close to the primes. We put the nontrivial zeros of the zeta function into $\omega$, and out came the primes! Another marvel!

There are also somewhat smaller spikes at prime powers, namely 4, 8, 9, 16, 25, 27, and 32.

The zeros of the zeta function know everything about primes. This also works in reverse: the primes know the zeros of the zeta function. To see this, consider the function
$$g_i(x)=\frac{\cos(x\cdot\ln(p_i))}{\sqrt{p_i}},$$
where $p_i$ is the $i$th prime. Again, we add these functions:
\[
\kappa(x)=\sum_{i=1}^{1000}g_i(x)=
\frac{\cos(x\cdot\ln(2))}{\sqrt2}
+\frac{\cos(x\cdot\ln(3))}{\sqrt3}
+\frac{\cos(x\cdot\ln(5))}{\sqrt5}
+\ldots+\frac{\cos(x\cdot\ln(7919))}{\sqrt{7919}}.
\]

\begin{figure}[htbp]
    \centering
    \includegraphics[width=0.8\textwidth]{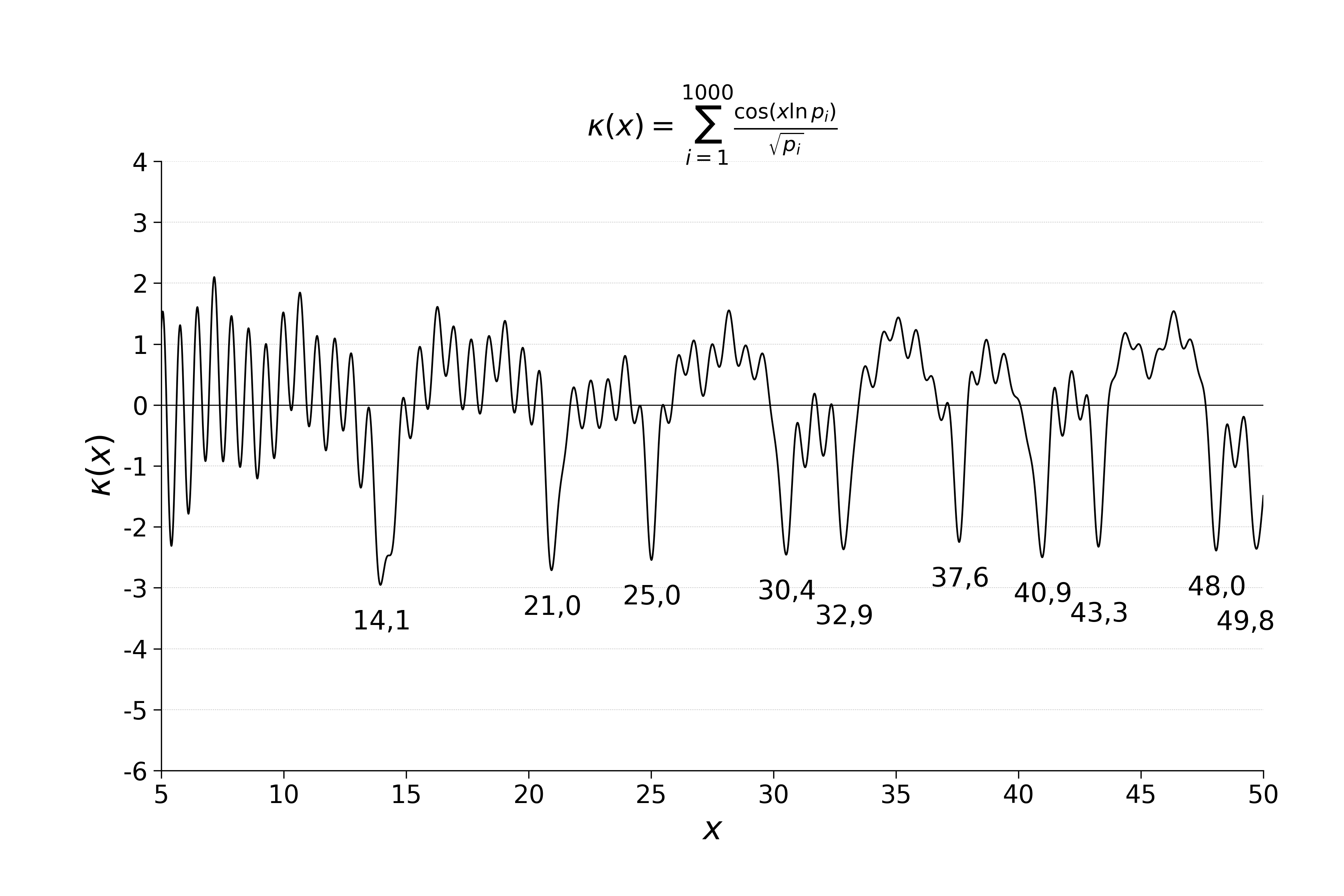}
    \caption{The function $\kappa(x)$.}
    \label{fig:kappa}
\end{figure}

Figure~\ref{fig:kappa} shows the function $\kappa(x)$.
The spikes lie very close to the imaginary parts of the nontrivial zeros of the zeta function. They do not occur at exactly those values because we have used only a finite sum for $\kappa$. A more precise mathematical treatment would show that we are investigating the zeros along the critical line.

The secret of the primes: they can be read off from the zeros of an analytic function.

\section{The Divisor-Sum Conjecture}

We define
\[\sigma(n)=\sum_{d|n}d.\]
Thus $\sigma(n)$ is the sum of all positive divisors of $n$. For example,
\[\sigma(10)=1+2+5+10=18.\]

It is known (and the proof is not difficult) that $\sigma$ is {\it multiplicative}: if $m,n$ are coprime natural numbers, then $\sigma(mn)=\sigma(m)\sigma(n)$. Thus, given the prime factorization $n=p_1^{k_1}p_2^{k_2}\ldots p_j^{k_j}$, we have
\[\sigma(n)=\sigma(p_1^{k_1})\sigma(p_2^{k_2})\ldots\sigma(p_j^{k_j}).\]
It is easy to see that if $p$ is prime, then $\sigma(p)=p+1$ and $\sigma(p^n)=1+p+p^2+\ldots+p^n$. 
Together with Formula (\ref{gf1}) we have
\[\sigma(n)=\prod_{i=1}^j\frac{1-p_i^{k_i+1}}{1-p_i}.\]
However, this formula allows us to compute $\sigma(n)$ only when we know the prime factorization of $n$.

We also need the function
\[h(n)=\sum_{i=1}^n\frac1i.\]
These are partial sums of the series considered in Theorem~\ref{srei}. The value $h(n)$ is called the {\it $n$th harmonic number}. In {\sc Julia}, the command

\verb+cumsum(1.0 ./ (1:n))+

gives us $h(n)$. First, \verb+1.0 ./ (1:n)+ produces the reciprocals $1,\frac12,\ldots,\frac1n$, and \verb+cumsum+ adds them cumulatively.

\begin{figure}[htbp]
    \centering
    \includegraphics[width=0.6\textwidth]{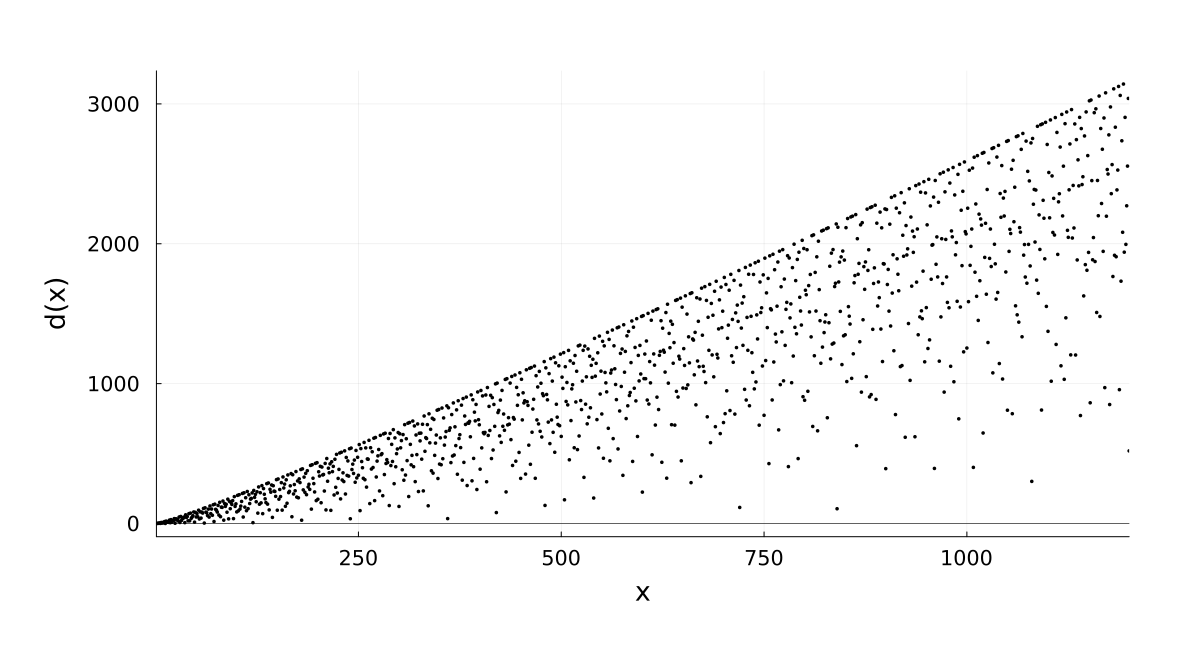}
    \caption{The divisor-sum conjecture.}
    \label{fig:TV}
\end{figure}

The divisor-sum conjecture states:
\begin{quote}
{\bf (DS):} For every natural number $n$, the following inequality holds:
\[\sigma(n)\le h(n)+e^{h(n)}\ln(h(n)).\]
\end{quote}
This conjecture can be formulated in entirely elementary terms.
In 2000, Lagarias surprisingly proved that it is equivalent to the Riemann hypothesis (see \cite{Bru01}). Using {\sc Julia}, we can plot the function
\[d(n)=h(n)+e^{h(n)}\ln(h(n))-\sigma(n).\]
Figure~\ref{fig:TV} shows $d(n)$. All values in this range are positive, as expected.
We can see how unevenly the values are distributed. They all appear to lie below a straight line.

\pagebreak

\section{Hardy's Function}

To conclude, let us experiment with {\sc Julia} once more.
Hardy's $Z$-function (see \cite{Iv16}) is a function $Z\colon\R\to\R$ whose zeros are exactly the imaginary parts of the nontrivial zeros of the zeta function that lie on the critical line (these are all the nontrivial zeros if the Riemann hypothesis holds). We have
\[
Z(t)=0\quad\Longleftrightarrow\quad\zeta\left(\frac12+it\right)=0.
\]
We cannot simply use $\zeta(1/2+ti)$ as our real-valued function, because away from its zeros it takes complex values. The $Z$-function is defined by
\[
Z(t)=e^{i\theta(t)}\,\zeta\!\left(\frac12+it\right),
\quad t\in\mathbb R,
\]
where
\[
\theta(t)=\operatorname{Im}\log\Gamma\!\left(\frac14+\frac{it}{2}\right)
-\frac t2\ln\pi
\]
is the Riemann--Siegel theta function.
Here $\log$ denotes the complex natural logarithm, and $\Gamma$ is the gamma function. I shall not go into the details here. What matters is that the factor $e^{i\theta(t)}$ represents a rotation (its modulus is 1) that rotates the complex function values onto the real axis, so that the resulting function takes only real values.

\begin{figure}[htbp]
    \centering
    \includegraphics[width=0.6\textwidth]{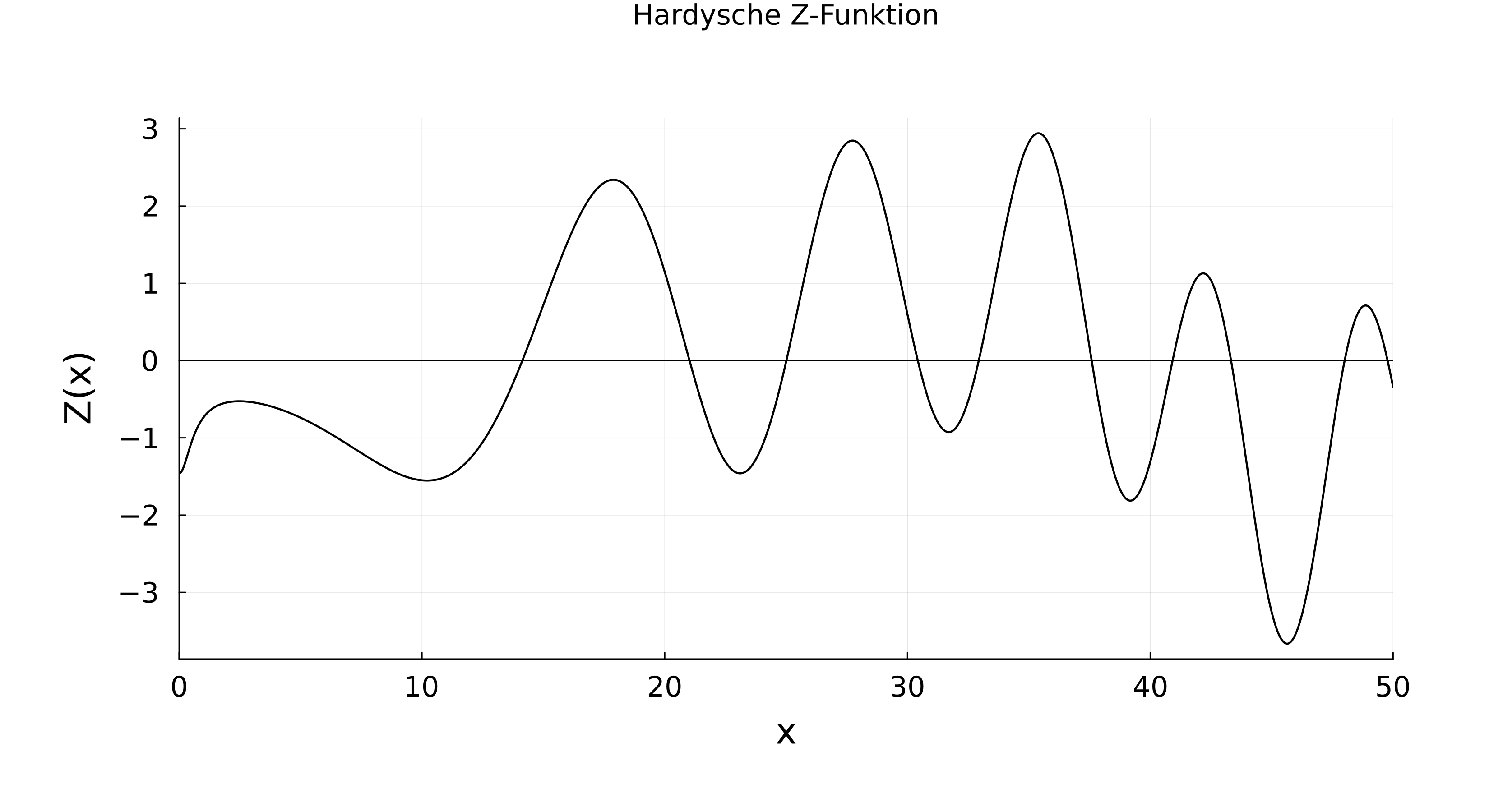}
    \caption{The $Z$-function from 0 to 50.}
    \label{fig:HZ50}
\end{figure}

The following program implements the $Z$-function in {\sc Julia}:
\begin{verbatim}
using SpecialFunctions

# Riemann-Siegel theta function
function theta(x::Real)
    return imag(loggamma(0.25 + 0.5im * x)) - (x / 2) * log(pi)
end

# Hardy's Z-function
function Z(x::Real)
    return real(exp(im * theta(x)) * zeta(0.5 + im * x))
end
\end{verbatim}

In Figure~\ref{fig:HZ50},
we have plotted part of the graph of the $Z$-function using the program above. We can see that its zeros correspond exactly to the first nontrivial zeros of the zeta function. The absolute values of the other function values agree exactly with the modulus of $\zeta(1/2+xi)$. In this sense, the $Z$-function visualizes the zeta function on the critical line. Looking at several plots over different ranges of the $x$-axis reveals how chaotic the function's behavior is.


AI was used in this article to generate some of the illustrations, to create the {\sc Julia} programs and the bibliography, and to proofread the text.

 Stephan Rosebrock\\
  Department of Mathematics\\
  Karlsruhe University of Education\\
  Bismarckstr. 10\\
  76133 Karlsruhe, Germany\\
  Email: rosebrock@ph-karlsruhe.de

\end{document}